\documentclass[11pt]{article}
\usepackage[dvips]{graphics}
\input{amssym.def}
\input{amssym.tex}

\title{\Large\bf The Baker-Akhiezer function and 
\\
factorization of the Chebotarev-Khrapkov matrix }
\author{Yuri A. Antipov \\
Department of Mathematics, Louisiana State University\\
Baton Rouge, LA 70803}

\date{}

\newcommand{\bfm}[1]{\mbox{\boldmath ${#1}$}}

\newcommand{\beqa}{\begin{eqnarray}}
\newcommand{\eeqa}[1]{\label{#1}\end{eqnarray}}
\newcommand{\bequ}{\begin{equation}}
\newcommand{\eequ}[1]{\label{#1}\end{equation}}

\newcommand{\Ga}{\alpha}
\newcommand{\Gb}{\beta}
\newcommand{\Gd}{\delta}

\newcommand{\Gve}{\varepsilon}
\newcommand{\Gf}{\phi}

\newcommand{\Gg}{\gamma}
\newcommand{\Gc}{\chi}

\newcommand{\Gk}{\kappa}

\newcommand{\Gl}{\lambda}

\newcommand{\Gt}{\theta}

\newcommand{\Gr}{\rho}

\newcommand{\Gs}{\sigma}

\newcommand{\Go}{\omega}
\newcommand{\Gx}{\xi}

\newcommand{\Gz}{\zeta}

\newcommand{\GF}{\Phi}

\newcommand{\GL}{\Lambda}

\newcommand{\GO}{\Omega}

\newcommand{\BGy}{\bfm\psi}

\newcommand{\BGF}{\bfm\Phi}

\newcommand{\CF}{{\cal F}}

\newcommand{\CL}{{\cal L}}

\newcommand{\CR}{{\cal R}}

\def\Ba{{\bf a}}
\def\Bb{{\bf b}}

\def\Bg{{\bf g}}

\def\Bn{{\bf n}}

\def\Bs{{\bf s}}

\def\BB{{\bf B}}

\newcommand{\beq}{\begin{equation}}
\newcommand{\eeq}{\end{equation}}
\newcommand{\barr}{\begin{eqnarray}}
\newcommand{\earr}{\end{eqnarray}}
\newcommand{\beqn}{\begin{equation*}}
\newcommand{\eeqn}{\end{equation*}}
\newcommand{\barrn}{\begin{eqnarray*}}
\newcommand{\earrn}{\end{eqnarray*}}

\newcommand{\fr}{\frac}

\newcommand{\diag}{\mbox{diag}}

\newcommand{\sgn}{\mathop{\rm sgn}\nolimits}

\newcommand{\ind}{\mathop{\rm ind}\nolimits}
\newcommand{\I}{\mathop{\rm Im}\nolimits}
\newcommand{\R}{\mathop{\rm Re}\nolimits}
\newcommand{\const}{\mbox{const}}

\begin{document}
\maketitle

\begin{abstract}

A new technique is proposed for the solution of the Riemann-Hilbert problem   with the Chebotarev-Khrapkov matrix coefficient
$G(t)=\Ga_1(t)I+\Ga_2(t)Q(t)$,
$\Ga_1(t), \Ga_2(t)\in H(L)$, $I=\diag\{1,1\}$, $Q(t)$ is a $2\times 2$ zero-trace polynomial matrix. This problem has numerous applications
in elasticity and diffraction theory. The main feature of the method is  the removal of essential singularities  of the solution
to the associated homogeneous scalar Riemann-Hilbert problem on the hyperelliptic surface of an algebraic function
by means of the Baker-Akhiezer function. The consequent application of this function for
the derivation of the general solution to the vector Riemann-Hilbert problem requires the finding of the $\Gr$ zeros of the Baker-Akhiezer function
($\Gr$ is the genus of the surface).
These zeros are recovered through the solution to the associated Jacobi problem of inversion of abelian integrals or, equivalently,
the determination of the  zeros of the associated  degree-$\Gr$ polynomial and solution of a certain linear algebraic system
of $\Gr$ equations.

\end{abstract}



\setcounter{equation}{0}

\section{Introduction}

Many problems of elasticity [\ref{khr}],  [\ref{moi2}], [\ref{ant5}], [\ref{ant1}], [\ref{ant2}],
electromagnetic diffraction [\ref{dan}], [\ref{hur}],  [\ref{buy}],  [\ref{lun}],  [\ref{ant7}],  [\ref{ant8}],  [\ref{ant3}],  [\ref{ant4}],   
 and acoustic diffraction  [\ref{jon}],  [\ref{raw}],  [\ref{ant6}]
 require the solution of the vector
 Riemann-Hilbert problem (RHP)
 of the theory of analytic functions [\ref{vek}]
$\BGF^+(t)=G(t)\BGF^-(t)+\Bg(t)$, $t\in L$,
where $L$ is either the whole real axis, or a finite segment, when the
 matrix $G(t)$ has the Chebotarev-Khrapkov (also known as Daniele-Khrapkov) 
 structure
 [\ref{che1}],   [\ref{khr}],   [\ref{dan}], 
 \beq
 G(t)=\Ga_1(t)I+\Ga_2(t)Q(t).
 \label{1.1}
\eeq
Here, $\Ga_1(t)$ and $\Ga_2(t)$ are H\"older functions on $L$, $I=\diag\{1,1\}$, and
$Q(t)$ is a $2\times 2$ zero-trace polynomial matrix. In the case $n=\deg f(z)\le 2$ ($\det Q(z)=h^2(z)f(z)$,
and $f(z)$ has simple zeros only)
the problem was solved in [\ref{khr}]. For a particular case of the matrix (\ref{1.1})  and when
 $n=4$, the exact solution was derived in [\ref{dan}]. For any finite $n$, 
 the vector problem is reduced [\ref{moi1}]  to a scalar RHP on a hyperelliptic surface
of genus $\Gr=[(n-1)/2]$.
A theory of the RHP on compact Riemann surfaces and a constructive procedure
for the solution of the associated Jacobi inversion problem was proposed in [\ref{zve}] 
(see also [\ref{zve2}]).
This technique was further developed and adjusted to specific needs of the RHPs
on hyperelliptic surfaces arising in elasticity [\ref{moi2}], [\ref{ant5}],
diffraction theory in  [\ref{ant6}], [\ref{ant7}], [\ref{ant8}]
 and for symmetric vector RHPs   in [\ref{ant3}], [\ref{ant4}]. 
The method for the vector RHP with the coefficient (\ref{1.1}) in the elliptic and hyperelliptic cases
first factorizes the coefficient of the associated scalar RHP using the Weierstrass
analogue of the Cauchy kernel. In general,
that solution has   an essential singularity at the infinite points of the surface due to unavoidable
poles of the Weierstrass kernel.
The next step of the procedure, the removal of the essential singularities,
leads to the classical problem of the inversion of abelian integrals and, eventually,
to the finding of the zeros of a certain degree-$\Gr$ polynomial.

The main goal of this paper was to develop a new factorization procedure for matrices of the form (\ref{1.1})
based on the use of the Baker-Akhiezer function.
The Baker-Akhiezer function plays an important role in the study of analytic properties of  eigenfunctions of ordinary differential operators with periodic coefficients
[\ref{dub1}], [\ref{its}], [\ref{dub3}], [\ref{kri}], [\ref{dub2}].
The representation of the 
Baker-Akhiezer function on a genus-$\Gr$ hyperelliptic surface $\CR$
\beq
\CF(P)=e^{\GO(P)}\fr{\Gt(u_1(P)-\Gs_1+V_1^\circ, \ldots, u_\Gr(P)-\Gs_\Gr+V_\Gr^\circ)} 
{\Gt(u_1(P)-\Gs_1, \ldots, u_\Gr(P)-\Gs_\Gr)}
\label{1.2}
\eeq
that we employ for the solution of the Wiener-Hopf matrix factorization problem was  first written by A. R. Its in context of the finite gap solutions of the KdV equation
[\ref{mat}]. Here, $P\in\CR$,
$\GO(P)$ is an abelian integral of the second kind with zero $A$-periods and a certain prescribed polynomial growth at the infinite point
of the surface $\CR$,  $\Gt$ is the theta Riemann function, $u_1,\ldots, u_\Gr$ form the canonical
basis of abelian integrals of the first kind,
$\Gs_j=k_j+u_j(P_1)+\ldots+u_j(P_\Gr)$,  $P_j$
are simple poles of the Baker-Akhiezer function,
$k_j$ are the Riemann constants associated with
the homology basis $\Ba_1$, $\Bb_1,\ldots, \Ba_\Gr, \Bb_\Gr$, and $V_j^\circ=(2\pi i)^{-1}\int_{\Bb_j}d\GO$,  $j=1,\ldots,\Gr$.

In section 2 we state the vector RHP in the real axis with the matrix coefficient (\ref{1.1}) and reduce it to a scalar RHP on a hyperelliptic
surface $\CR$ of the algebraic function $w^2=f(z)$. We derive a particular solution, $\Gc_0(z,w)$, to the scalar RHP in section 3. This solution
satisfies the boundary condition but has inadmissible essential singularities at the two infinite points $\infty_1$ and $\infty_2$ of the surface.
In section  4 we construct the Baker-Akhiezer function (\ref{1.2}) of the surface $\CR$. This function
is associated with an
abelian integral of the second type with zero-$A$-periods  used to remove the essential singularities and two Riemann $\Gt$-functions which serve to make the solution
continuous through the $B$-cross-sections.  We find the Wiener-Hopf factors in terms of the functions $\Gc_0(z,w)$ and $\CF(P)$ and the general solution to the vector RHP
in section 5.

\setcounter{equation}{0}

\section{Scalar RHP on a Riemann surface associated  with  the Chebotarev-Khrapkov matrix}

Motivated by numerous applications in acoustics, electromagnetic theory, fluid mechanics and elasticity we assume that the Riemann-Hilbert contour,
 $L$, is the whole real axis which splits the plane of a complex variable $z$ into two half-planes, $D^+:\I z>0$ and $D^-:\I z<0$. Let 
$G(t)$ be a $2\times 2$ matrix which is nonsingular in $L$ and whose structure  is
\beq
G(t)=\left(\begin{array}{cc}
\Ga_1(t)+\Ga_2(t)l_0(t) & \Ga_2(t)l_1(t) \\
 \Ga_2(t)l_2(t) & \Ga_1(t)-\Ga_2(t)l_0(t)\\
\end{array}\right),
\label{2.1}
\eeq
where $\Ga_1(t), \Ga_2(t)l_j(t)\in \hat H(L)$, $j=0,1,2$, $l_0(t)$, $l_1(t)$ and $l_2(t)$ are polynomials, and $\hat H(L)$ is the class of all H\"older functions $\Ga(t)$ in
any finite interval in $L$ which tend to a definite limit $\Ga(\infty)$ as $t\to\pm\infty$.
For large $t$, they satisfy the condition $|\Ga(t)-\Ga(\infty)|<C|t|^{-\mu}$, $\mu>0$, $C>0$.
Without loss of generality assume that $\det G(\infty)=1$.
Let 
$\Bg(t)$ be an order-2 $\hat H$-vector-function on $L$ such that $\Bg(\infty)$ is the zero-vector.
Consider the following RHP.

{\it Given $G(t)$ and $\Bg(t)$ find two vectors, $\BGF^+(z)$ and $\BGF^-(z)$, analytic in the domains $D^+$ and $D^-$, respectively,  bounded at infinity,
$\hat H$-continuous up to the contour  $L$
and satisfying the boundary condition
\beq
\BGF^+(t)=G(t)\BGF^-(t)+\Bg(t), \quad t\in L.
\label{2.2} 
\eeq
}

Denote $l_0^2(z)+l_1(z) l_2(z)=h^2(z)f(z)$ and $f(z)=z^n+\Gve_1 z^{n-1}+\ldots+\Gve_n$.
All zeros,   $r_1, r_2,\ldots,r_n$, of $f(z)$ are simple, while the zeros of the polynomial
$h(z)$,  $p_1,p_2,\ldots,p_l$, have multiplicity $m_1, m_2, \ldots, m_l$, respectively, and $m_1+m_2+\ldots+m_l=N$. 
Some or all zeros of the polynomial $l_0^2(z)+l_1(z) l_2(z)$ may
have an odd multiplicity $2m_i+1\ge 3$. In this case the $i$-th zero is counted as a simple zero of $f(z)$ and an order-$m_i$ zero of the polynomial $h(z)$.
Assume  that none of the zeros of $f(z)$ and $h(z)$ falls in the contour $L$ (we refer to  [\ref{ant1}] otherwise).
In addition, we assume that $n$ is even,  $n=2\Gr+2$ (this is true for all known applications
of the problem (\ref{2.2}) with the matrix coefficient (\ref{2.1}) to elasticity and diffraction theory). This implies $\deg[l_0^2(z)+l_1(z) l_2(z)]=2N+2\Gr+2$. Denote $\deg l_j(z)=\Gd_j$, $j=0,1,2$, and for simplicity, 
accept that
 $0\le \Gd_0\le N+\Gr+1$ and $0\le \Gd_j\le 2N+2\Gr+2$, $j=1,2$ ($\Gd_1+\Gd_2\le 2N+2\Gr+2$).

Choose a single branch of $f^{1/2}(z)$ in the plane cut along simple smooth disjoint arcs
$\Gg_1=r_1r_2$, $\Gg_2=r_3r_4, \ldots$,  $\Gg_{\Gr+1}=r_{2\Gr+1}r_{2\Gr+2}$ such that  $f^{1/2}(z)\sim  z^{\Gr+1}$, $z\to\infty$.
 The functions 
 \beq
 \Gl_1(t)=\Ga_1(t)+\Ga_2(t)h(t)\sqrt{f(t)},\quad 
 \Gl_2(t)=\Ga_1(t)-\Ga_2(t)h(t)\sqrt{f(t)}
 \label{2.2'}
 \eeq
  are the eigenvalues
of the matrix $G(t)$, and their product $\Ga_1^2(t)-\Ga_2^2(t)h^2(t)f(t)$ is the determinant of $G(t)$.
To pursue the Wiener-Hopf factorization of $G(t)$, we split it as 
 \beq
 G(t)=T(t)\GL(t)[T(t)]^{-1},
 \label{2.5}
 \eeq
 where $\GL(t)=\diag\{\Gl_1(t),\Gl_2(t)\}$,
 \beq
T(z)=\left(\begin{array}{cc}
1 \; & 1 \\
 -\fr{l_0(z)-h(z)\sqrt{f(z)}}{l_1(z)} \;&  -\fr{l_0(z)+h(z)\sqrt{f(z)}}{l_1(z)} \\
\end{array}\right),
\label{2.6}
\eeq
and reduce the problem of matrix factorization to a scalar RHP on a Riemann surface  [\ref{moi1}].  
 First we introduce two new vectors, $\BGy(z)=(\psi_1(z), \psi(z))$ and $\Bg^\circ(t)=(g_1^\circ(t),g_2^\circ(t))$,
 \beq
 \BGy(z)=[T(z)]^{-1}\BGF(z), \quad \Bg^\circ(t)=[T(t)]^{-1}\Bg(t),
 \label{2.7}
 \eeq
where
 \beq
[T(z)]^{-1}=\left(\begin{array}{cc}
\fr{l_0(z)}{2h(z)\sqrt{f(z)}}+\fr12 \; & \fr{l_1(z)}{2h(z)\sqrt{f(z)}} \\
-\fr{l_0(z)}{2h(z)\sqrt{f(z)}}+\fr12 \; & -\fr{l_1(z)}{2h(z)\sqrt{f(z)}} \\
\end{array}\right).
\label{2.7'}
\eeq
The components of the vector $\BGy(z)$  are expressed 
through the components of the vector $\BGF(z)$ as
$$
\psi_1(z)=\fr12\left[1+\fr{l_0(z)}{h(z)\sqrt{f(z)}}\right]\GF_1(z)+\fr{l_1(z)}{2h(z)\sqrt{f(z)}}\GF_2(z),
$$
\beq
\psi_2(z)=\fr12\left[1-\fr{l_0(z)}{h(z)\sqrt{f(z)}}\right]\GF_1(z)-\fr{l_1(z)}{2h(z)\sqrt{f(z)}}\GF_2(z).
\label{2.7''}
\eeq
Similar formulas can be written for the components of the vectors $\Bg^\circ(t)$ and $\Bg(t)$. 
The new functions $\psi_1(z)$ and $\psi_2(z)$ may grow at infinity if $\Gd_1> N+\Gr+1$.
Let $\Gd=\max\{0,\Gd_1-N-\Gr-1\}$. Then since the functions $\GF_1(z)$ and $\GF_2(z)$ are bounded as $z\to\infty$, we have
$|\psi_j(z)|<c_j|z|^{\Gd}$, $z\to\infty$, $c_j=\const$,  $j=1,2$.

Due to continuity of the vector $\BGF(z)$ through the branch cuts $\Gg_j$ ($j=1,2,\ldots,\Gr+1$),
we have $T^+(t)\BGy^+(t)=T^-(t)\BGy^-(t)$, $t\in\Gg_j$. This implies that the components
of the vector $\BGy(z)$ satisfy the following Riemann-Hilbert boundary conditions:
$$
\psi_1^+(t)=\psi_2^-(t), \quad  \psi_2^+(t)=\psi_1^-(t), \quad t\in\Gg_j, \quad  j=1,2,\ldots,\Gr+1,
 $$
 \beq
 \psi_j^+(t)=\Gl_j(t)\psi_j^-(t)+g_j^\circ(t), \quad t\in L, \quad j=1,2,
 \label{2.8}
 \eeq
and may have poles $p_1,p_2,\ldots,p_l$ of multiplicity $m_1,m_2$, $\ldots,m_l$ 
at the zeros of the polynomial $h(z)$.

We wish to reformulate (\ref{2.8}) as a scalar RHP on a Riemann surface. Let $\CR$ be the two-sheeted
 Riemann surface of the algebraic function $w^2=f(z)$ formed by gluing two copies, $\Bbb C_1$ and $\Bbb C_2$,
 of the extended complex plane $\Bbb C\cup \infty$ along the cuts $\Gg_j$ ($j=1,2,\ldots,\Gr+1$) such that
 \beq
 w=\left\{\begin{array}{cc}
\sqrt{f(z)}, & z\in{\Bbb C_1}, \\
-\sqrt{f(z)}, & z\in{\Bbb C_2}, \\
\end{array}\right.
\label{2.9}
\eeq
 is a single-valued function on the surface $\CR$. Here, $\sqrt{f(z)}$ is the branch chosen before. 
 Let $\Ba_j$, $\Bb_j$ ($j=1,2,\ldots,\Gr$) be a homology basis of the genus-$\Gr$ surface $\CR$ 
 (Figure 1).
 \begin{figure}[t]
\centerline{
\scalebox{0.7}{\includegraphics{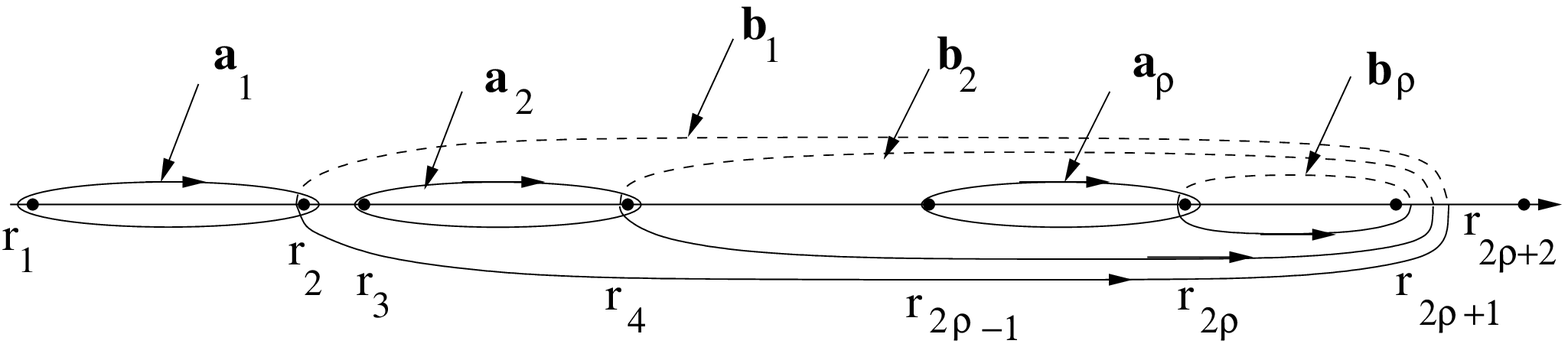}}}
\caption{The canonical cross-sections $\Ba_j$ and $\Bb_j$, $j=1,\ldots,\Gr$.}
\label{fig1}
\end{figure} 
 Denote $\CL=L_1\cup L_2$ the contour on the surface $\CR$ with $L_j\subset{\Bbb C}_j$ ($j=1,2$)
being two copies of the contour $L$. With each pair of the functions $(\psi_1,\psi_2)$, $(\Gl_1,\Gl_2)$ and 
 $(g_1^\circ, g_2^\circ)$ we associate the following functions on the surface $\CR$:
 $$
 \Psi(z,w)=\psi_j(z), \quad (z,w)\in{\Bbb C}_j,
 $$
 \beq
 \Gl(t,\Gx)=\Gl_j(t), \quad g^*(t,\xi)=g_j^\circ(t), \quad (t,\Gx)\in L_j, \quad j=1,2, \quad \Gx=w(t).
 \label{2.10}
 \eeq
The function $\Psi(z,w)$ may have simple poles at the branch points of the surface $\CR$, $r_1,r_2, \dots, r_{2\Gr+2}$ (recall  [\ref{spr}] that a branch
point $r_j$ of the Riemann surface $\CR$ is called an order-$l_j$ pole of the function $\Psi(z,w)$ if
$\Psi(z,w)\sim A_j\Gz^{-l_j} $, $\Gz\to 0$, $A_j=\;$const, and $\Gz = (z-r_j)^{1/2}$ is a local uniformizing
parameter of the point $r_j$).
We also assert that the function $\Psi(z,w)$ is continuous through the contours $\Gg_j$ ($j=1,2,\ldots,\Gr+1$), and therefore
the vector RHP (\ref{2.2}) on the plane is equivalent to the following scalar RHP on the surface $\CR$.

{\it Find a piece-wise analytic function $\Psi(z,w)$ with the discontinuity contour $\CL$,
$\hat H$-continuous up to the contour
$\CL$, satisfying the boundary condition
\beq
\Psi^+(t,\xi)=\Gl(t,\xi)\Psi^-(t,\xi)+g^*(t,\xi), \quad (t,\xi)\in\CL,
\label{2.11}
 \eeq
and having  poles $p_1,p_2,\ldots,p_l$ of multiplicity $m_1,m_2$, $\ldots,m_l$
in both sheets of the surface $\CR$ and simple poles at the branch points $r_1,r_2, \dots, r_{2\Gr+2}$. In neighborhoods of 
 the two infinite points $\infty_j=(\infty, (-1)^{j-1}\infty)$  of the surface $\CR$ the function $\Psi(z,w)$ satisfies the inequality
$ |\Psi(z,w)|<c_j |z|^\Gd$, $c_j=\const$, $j=1,2$.}

\setcounter{equation}{0}

\section{Solution with an essential singularity at the infinite points}

We begin with factorization of the function $\Gl(t,\xi)$. For an analogue of the Cauchy kernel we choose the Weierstrass kernel
\beq
dW=\fr{w+\Gx}{2\Gx}\fr{dt}{t-z}
\label{3.1}
\eeq
and analyze the integral
$$
\fr{1}{2\pi i}\int_\CL\log\Gl(t,\xi)dW=\fr{1}{4\pi i}\int_L\left[\log\Gl_1(t)+\log\Gl_2(t)\right]\fr{dt}{t-z}
$$
\beq
+\fr{w}{4\pi i}\int_L\left[\log\Gl_1(t)-\log\Gl_2(t)\right]\fr{dt}{\sqrt{f(t)}(t-z)}.
\label{3.2}
\eeq
Pick a point on $L$, $z_0$, and treat it as the starting point, $z_0^+$, of the contour $L$ (it is convenient to take $z_0=0$).  
Let
\beq
\Gk_j=\ind\Gl_j(t)=\fr{1}{2\pi}[\arg\Gl_j(t)]|_{L},
\label{3.3}
\eeq
where $\ind \Gl_j(t)$ is the index of the function $\Gl_j(t)$, and
$[\arg\Gl_j(t)]|_L$ is the increment of $\arg\Gl_j(t)$ as $t$ traverses the contour $L$ in the positive direction with $z_0^+$ being the starting point.
Because of the continuity of the functions $\Gl_1(t)$ and $\Gl_2(t)$ in the contour $L$ both numbers, $\Gk_1$  and $\Gk_2$, are  integers.
Fix branches of the logarithmic functions $\log\Gl_1(t)$
and $\log\Gl_2(t)$ by the condition
$\arg\Gl_j(z_0^+)=\Gf_j$, $ 0\le\Gf_j<2\pi.$
Then at the terminal point $z_0$ of the contour $L$ (to distinguish the  terminal and starting points, we denote the former point as $z_0^-$),
$\arg\Gl_j(z_0^-)=\Gf_j+2\pi\Gk_j$. Analysis of the singular integrals in the right-hand side (\ref{3.2}) implies
$$
\fr{1}{4\pi i}\int_L[\log\Gl_1(t)+\log\Gl_2(t)]\fr{dt}{t-z}\sim\fr{\Gk_1+\Gk_2}{2}\log(z-z_0),\quad z\to z_0,
$$
\beq
\fr{w}{4\pi i}\int_L[\log\Gl_1(t)-\log\Gl_2(t)]\fr{dt}{\sqrt{f(t)}(t-z)}\sim\fr{\Gk_1-\Gk_2}{2}(-1)^{j-1}\log(z-z_0),\quad z\to z_0, \quad (z,w)\in{\Bbb C}_j.
\label{3.6}
\eeq
Consequently, the integral in the left-hand side
(\ref{3.2}) has a logarithmic singularity at the point $(z_0,w(z_0))\in\CL$ in both sheets of the surface
\beq
\fr{1}{2\pi i}\int_\CL\log\Gl(t,\xi)dW\sim \Gk_j\log(z-z_0), \quad  z\to z_0, \quad (z,w)\in{\Bbb C}_j,  \quad j=1,2.
\label{3.7}
\eeq
It is an easy matter to move the singularity from the contour $\CL$ to the surface $\CR\setminus\CL$
by adding the extra term
\beq
I(z,w)=\sum_{m=1}^2\sgn\Gk_m\sum_{j=1}^{|\Gk_m|}
\int_{q_{m 0}}^{q_{m j}}dW.
\label{3.9}
\eeq
 Here, $q_{m0}q_{mj}\subset{\Bbb C_m}$ are smooth simple contours which do not intersect the contours $L_m$,
 $q_{m j}=(z_{mj}, (-1)^{m-1}\sqrt{f(z_{mj})})\in{\Bbb C}_m\setminus L_m$,
$j=1,2,\ldots,|\Gk_m|$, are arbitrary fixed points, $z_{mj}$ are their affixes, and $q_{m0}=(z_0,(-1)^{m-1}\sqrt{f(z_0)})\in L_m$, $m=1,2$.  
The function $\exp\{I(z,w)\}$ is continuous through the contour $\CL$ except for the points $q_{10}$ and $q_{20}$ at which the integral $I(z,w)$
has logarithmic singularities. In  addition, the integral $I(z,w)$ has logarithmic singularities at the internal points 
$q_{m j}$,
$$
I(z,w)\sim -\Gk_m\log(z-z_0), \quad  z\to z_0, \quad (z,w)\in{\Bbb C}_m,  
$$
\beq
I(z,w)\sim\sgn\Gk_m\log(z-z_{mj}), \quad  z\to z_{mj}, \quad (z,w)\in{\Bbb C}_m,  \quad j=1,\ldots,|\Gk_m|, \quad m=1,2.
\label{3.10}
\eeq
At the same time, the sum of the integral (\ref{3.2}) and  (\ref{3.9})  does not have the singularity at the points
$(z_0,\pm\sqrt{f(z_0)})$. 
Now, to factorize the function $\Gl(t,\xi)$, we use
the function
\beq
\Gc_0(z,w)=\exp\left\{\fr{1}{2\pi i}\int_\CL\log\Gl(t,\xi)dW+\sum_{m=1}^2\sgn\Gk_m\sum_{j=1}^{|\Gk_m|}
\int_{q_{m 0}}^{q_{m j}}dW\right\}.
\label{3.11}
\eeq
The function $\Gc_0(z,w)$ satisfies the homogeneous boundary condition
\beq
\Gc_0^+(t,\Gx)=\Gl(t,\xi)\Gc_0^-(t,\Gx), \quad (t,\xi)\in \CL,
\label{3.13}
\eeq
does not have singularities  at the points  $q_{10}$ and $q_{20}$,
but has inadmissible essential singularities at the points $\infty_1$ and $\infty_2$. 
Also,  it has simple zeros $z_{mj}$
on the sheet ${\Bbb C}_m$ if $\Gk_m>0$ and simple poles  $z_{mj}$ on ${\Bbb C}_m$ if  $\Gk_m<0$ ($j=1,\ldots,|\Gk_m|$, $m=1,2$).

\setcounter{equation}{0}

\section{Baker-Akhiezer function}

Our aim is to quench the essential singularities at the infinite points of the function $\Gc_0(z,w)$ by employing  the Baker-Akhiezer function,  $\CF(z,w)$,  on the genus-$\Gr$ 
surface $\CR$ associated with the function $\Gc_0(z,w)$. The function $\CF(z,w)$
has to satisfy the following two conditions: 

(i) it is meromorphic everywhere on $\CR$ except at the points $\infty_1$ and $\infty_2$,

(ii) the function $\Gc_0(z,w)\CF(z,w)$ is bounded at the points $\infty_1$ and $\infty_2$.

Setting
\beq
\Gc_0(z,w)=e^{\Gb_0(z)+w\Gb_1(z)},
\label{4.00}
\eeq
where 
$$
\Gb_0(z)=\fr{1}{4\pi i}\int_L\left[\log\Gl_1(t)+\log\Gl_2(t)\right]\fr{dt}{t-z}+\fr12\sum_{m=1}^2\sgn\Gk_m\sum_{j=1}^{|\Gk_m|}
\int_{z_0}^{z_{mj}}\fr{dt}{t-z},
$$
\beq
\Gb_1(z)=\fr{1}{4\pi i}\int_L\left[\log\Gl_1(t)-\log\Gl_2(t)\right]\fr{dt}{\sqrt{f(t)}(t-z)}+\fr12\sum_{m=1}^2\sgn\Gk_m\sum_{j=1}^{|\Gk_m|}
\int_{q_{m0}}^{q_{mj}}\fr{dt}{\Gx(t-z)},
\label{4.0}
\eeq
we study the behavior of the function $\Gc_0(z,w)$ at the infinite points.
For  the branch $\sqrt{f(z)}$ chosen we have
 \beq
 \sqrt{f(z)}=\sqrt{\prod_{j=1}^{2\Gr+2}(z-r_j)}=z^{\Gr+1}\sum_{m=0}^\infty c_m z^{-m},
\label{4.1}
\eeq
Here,
$$
c_0=1,\quad 
c_1=\fr{(-1/2)_1}{1!}\sum_{j=1}^{2\Gr+2} r_j,
$$$$
c_2=\fr{(-1/2)_2}{2!}\sum_{j=1}^{2\Gr+2}r_j^2+\fr{[(-1/2)_1]^2}{(1!)^2}\sum_{j=1}^{2\Gr+2}r_j\sum_{m=1,m\ne j}^{2\Gr+2}r_m,
$$
\beq
c_3=\fr{(-1/2)_3}{3!}\sum_{j=1}^{2\Gr+2}r_j^3+\fr{(-1/2)_1(-1/2)_2}{1!2!}\sum_{j=1}^{2\Gr+2}r_j^2\sum_{m=1,m\ne j}^{2\Gr+2}r_m, \ldots,
\label{4.2}
\eeq
where $(a)_m=a(a+1)\ldots(a+m-1)$ is the factorial symbol. By virtue of (\ref{4.0})
\beq
\Gb_1(z)=\sum_{j=0}^\infty\fr{\tilde c_j}{z^{j+1}},
\label{4.3}
\eeq
where
\beq
\tilde c_j=-\fr{1}{4\pi i}\int_L[\log\Gl_1(t)-\log\Gl_2(t)]
\fr{t^jdt}{\sqrt{f(t)}}
-\fr12\sum_{m=1}^2\sgn\Gk_m
\sum_{j=1}^{|\Gk_m|}
\int_{q_{m0}}^{q_{mj}}
\fr{t^jdt}{\Gx},
\label{4.4}
\eeq
and therefore, as $z\to\infty$,
\beq
\sqrt{f(z)}\Gb_1(z)=z^{\Gr}\sum_{m=0}^\infty\fr{d_m}{z^m}, \quad d_m=\sum_{k=0}^m c_k\tilde c_{m-k}.
\label{4.5}
\eeq
This brings us to the expansion of the function $\Gc_0(z,w)$ at the infinite points
\beq
\Gc_0(z,w)=\exp\{(-1)^{j-1}M(z)+O(1)\}, \quad (z,w)\to\infty_j, \quad (z,w)\in{\Bbb C_j}, \quad j=1,2,
\label{4.6}
\eeq
where
\beq
M(z)=d_0z^\Gr+d_1 z^{\Gr-1}+\ldots+d_{\Gr-1}z.
\label{4.7}
\eeq

Our next step is to construct a special abelian integral of the second kind, 
\beq
\GO(P)=\int_{P_0}^P d\GO, \quad P_0=(r_{2\Gr+2},0), \quad P=(z,w).
\label{4.8}
\eeq
Determine $\GO(P)$ by the following properties:

(a) $\GO(P)\sim (-1)^j  M(z)$, $P\to\infty_j\in{\Bbb C}_j$, $j=1,2$,

(b) $ \int_{\Ba_j}d\GO=0, \quad j=1,2,\ldots,\Gr$. 

We seek the abelian differential $d\GO$ in the form
\beq
d\GO=\fr{e_0z^{2\Gr}+e_1 z^{2\Gr-1}+\ldots +e_{2\Gr}}{w}dz,
\label{4.9}
\eeq
where the coefficients $e_j$ are to be determined. We wish to exploit this formula in order to study the behavior of the integral
$\GO(P)$  at the infinite points. Because of (\ref{4.1}) we have
\beq
d\GO=(-1)^{j-1}\left(\tilde e_0 z^{\Gr-1}+\tilde e_1 z^{\Gr-2}+\ldots+\tilde e_{\Gr-1}+\fr{\tilde e_\Gr}{z}+\ldots\right)dz,
\label{4.10}
\eeq
where $\tilde e_m$  are defined recursively by
\beq
\tilde e_m=-\sum_{k=1}^m \tilde e_{m-k} c_k +e_m, \quad m=0,1,\ldots,\Gr.
\label{4.11}
\eeq
By integrating (\ref{4.10}) we determine the asymptotic expansion  of the abelian integral $\GO(P)$
$$
\GO(P)=(-1)^{j-1}\left(\fr{\tilde e_0 z^{\Gr}}{\Gr}+\fr{\tilde e_1 z^{\Gr-1}}{\Gr-1}+\ldots+\tilde e_{\Gr-1}z
+\tilde e_\Gr \log z-\fr{\tilde e_{\Gr+1}}{z}-\ldots\right)+K_j,
$$
\beq
 P\to\infty_j\in{\Bbb C}_j, \quad j=1,2,
\label{4.12}
\eeq
where $K_1$  and $K_2$ are constants. 
On satisfying the property (a) of the integral $\GO(P)$ we find the coefficients $\tilde e_0$, \ldots,$\tilde e_\Gr$
\beq
\tilde e_0=-\Gr d_0, \quad \tilde e_1=-(\Gr-1)d_1, \quad \tilde e_2=-(\Gr-2)d_2, \; \ldots,\;
\tilde e_{\Gr-1}=-d_{\Gr-1}, \quad \tilde e_\Gr=0.
\label{4.14}
\eeq
Due to (\ref{4.11}) we can express the coefficients $e_m$ ($m=0,1,\ldots,\Gr$) through $\tilde e_m$
\beq
e_m=\tilde e_m+\sum_{k=1}^m \tilde e_{m-k} c_k, \quad m=0,1,\ldots,\Gr.
\label{4.16}
\eeq
The remaining coefficients $e_{\Gr+1}, e_{\Gr+2}, \ldots, e_{2\Gr}$ in the representation (\ref{4.9}) of the abelian differential are fixed by solving
the system of $\Gr$ linear algebraic equations
\beq
\sum_{m=\Gr+1}^{2\Gr} U_{jm} e_m=\hat d_j, \quad j=1,2,\ldots,\Gr,
\label{4.17}
\eeq
which follows from the property (b) of the integral $\GO(P)$. Here, 
\beq
\hat d_j=-\sum_{m=0}^\Gr U_{jm} e_m, \quad U_{jm}=\int_{\Ba_j}\fr{z^{2\Gr-m}}{w(z)}dz.
\label{4.18}
\eeq
This completes the construction of the abelian integral $\GO(P)$. 

It becomes evident that the product $\Gc_0(z,w)\exp\{\GO(P)\}$ is bounded as $P\to\infty_j\in{\Bbb C}_j$, $j=1,2$.  This function is continuous through
the cross-sections $\Ba_j$ of the surface $\CR$ because of the zero $A$-periods and discontinuous through the cross-sections
$\Bb_j$ ($j=1,2,\ldots,\Gr$) due to the non-zero $B$-periods of the integral $\GO(P)$. Our efforts will now be directed towards annihilating 
the jumps $\exp\{V_m\}$, 
\beq
V_m=\int_{\Bb_m}d\GO, \quad m=1,2,\ldots,\Gr,
\label{4.19}
\eeq
 of the function $\exp\{\GO(P)\}$
through the cross-sections $\Bb_m$, $m=1,2,\ldots,\Gr$.

Let $d\Go_j$ ($j=1,2,\ldots,\Gr$) be the canonical basis of Abelian differentials of the first kind
\beq
d\Go_j=\fr{c_j^{(1)}z^{\Gr-1}+c_j^{(2)}z^{\Gr-2}+\ldots+c_j^{(\Gr)}}{w}dz,
\label{4.20}
\eeq
where the constants $c_j^{(k)}$ ($k,j=1,2,\ldots,\Gr$) are chosen such that
\beq
\int_{\Ba_k}d\Go_j =\Gd_{jk}.
\label{4.21}
\eeq
Denote the $B-$periods of the basis $d\Go_j$ by
\beq
B_{jk}=\int_{\Bb_k}d\Go_j.
\label{4.22} 
\eeq
The matrix $B=(B_{jk})$ ($j,k=1,2,\ldots,\Gr$) is symmetric and $\I B$ is positive definite.
The principal tool we shall use to suppress the discontinuities of $\exp\{\GO(P)\}$
is the Riemann $\Gt$-function
\beq
\Gt(\Bs(P))=\Gt(s_1(P),s_2(P),\ldots,s_\Gr(P))
\label{4.23}
\eeq
defined by
\beq
\Gt(\Bs(P))=\sum_{n_1,\ldots,n_\Gr=-\infty}^\infty\exp\left\{
\sum_{j=1}^\Gr\sum_{k=1}^\Gr
B_{jk}n_j n_k+2\pi i\sum_{j=1}^\Gr n_j s_j(P)\right\}.
\label{4.24}
\eeq
Because of the positive definiteness of the matrix $\I B$ the series converges for all $\Bs(P)$.
The $\Gt$-function has periods $\Bn=(n_1,n_2,\ldots,n_\Gr)$, $n_j$ are integers, and
quasiperiods $\BB_j=(B_{j1}, B_{j2}, \ldots, B_{j\Gr})$, $j=1,2,\ldots,\Gr$,
$$
\Gt(s_1+n_1, \ldots, s_\Gr+n_\Gr)=\Gt(s_1, \ldots, s_\Gr),
$$
\beq
\Gt(s_1+B_{j1}, \ldots, s_\Gr+B_{j\Gr})=
\exp\{-\pi i B_{jj}-2\pi i s_j\} \Gt(s_1, \ldots,s_\Gr).
\label{4.25}
\eeq
Introduce next the function
\beq
\CF_0(P)=\fr{\Gt(u_1(P)-\Gs_1+V_1^\circ, \ldots, u_\Gr(P)-\Gs_\Gr+V_\Gr^\circ)} 
{\Gt(u_1(P)-\Gs_1, \ldots, u_\Gr(P)-\Gs_\Gr)}.
\label{4.26}
\eeq
Here,  $V_j^\circ=(2\pi i)^{-1}V_j$ and $u_j(P)$ are the integrals
\beq
u_j(P)=\int_{P_0}^P d\Go_j, \quad j=1,2,\ldots,\Gr,
\label{4.27}
\eeq
which form the canonical basis of  abelian integrals of the first kind. It is convenient to choose $P_0$ as the branch point $r_{2\Gr+2}$.
The numbers $\Gs_j$ are chosen to be
\beq
\Gs_j=\sum_{m=1}^\Gr u_j(P_m)+k_j,\quad j=1,\ldots,\Gr,
\label{4.28}
\eeq
where $P_m$ ($m=1,2,\dots,\Gr$) are some arbitrary distinct fixed 
points on $\CR$ say, on ${\Bbb C}_1$, $P_m=(\Gz_m,\sqrt{f(\Gz_m)})$, such that the $\Gt$-functions in (\ref{4.26}) are not identically zero.
The parameters $k_j$ in (\ref{4.28}) are the Riemann constants which, for the hyperelliptic surface $\CR$
and for the homology basis  chosen, can be taken as (see for example [\ref{ant8}])
\beq
k_j=-\fr{j}{2}+\fr12\sum_{m=1}^\Gr B_{jm}.
\label{4.29} 
\eeq
The function $\CF_0(P)$ has $\Gr$ simple 
poles $P_1$, $P_2,\ldots,, P_\Gr$ [\ref{che2}, p. 303] lying in the first sheet and $\Gr$ simple zeros which may lie on either sheet of the surface. Call these zeros
$Q_j=(t_j,w_j)$,  where $w_j=\sqrt{f(t_j)}$ if $Q_j\in{\Bbb C}_1$ and  $w_j=-\sqrt{f(t_j)}$ if $Q_j\in{\Bbb C}_2$, $j=1,2,\ldots,\Gr$. The 
position of these zeros is unknown {\it a priori}, and without loss of generality these zeros are assumed to be simple.
The function $\CF_0(P)$  is continuous through the cross-sections $\Ba_j$ and discontinuous through the cross-sections
$\Bb_j$, $j=1,\ldots,\Gr$.  Due to (\ref{4.25}) its jumps are $\exp\{-V_j\}$. This implies that the function
$\CF(z,w)=\exp\{\GO(P)\}\CF_0(P)$ is meromorphic  on $\CR$ (it is continuous through the loops $\Bb_j$). The set of
singularities of the function $\CF(z,w)$ comprises the two infinite points $\infty_1$ and $\infty_2$
and $\Gr$ simple poles $P_m=(\Gz_m,\sqrt{f(\Gz_m)})\in{\Bbb C}_1$,  $m=1,2,\dots,\Gr$.
Therefore
\beq
\CF(P)=e^{\GO(P)}\fr{\Gt(u_1(P)-\Gs_1+V_1^\circ, \ldots, u_\Gr(P)-\Gs_\Gr+V_\Gr^\circ)} 
{\Gt(u_1(P)-\Gs_1, \ldots, u_\Gr(P)-\Gs_\Gr)}
\label{4.30}
\eeq
is the Baker-Akhiezer function of the surface $\CR$ with the homology basis $\Ba_j$, $\Bb_j$
($j=1,\ldots,\Gr$) associated with the abelian integral $\GO(P)$ and the poles $P_1,\ldots,P_\Gr$.

\setcounter{equation}{0}
\section{Vector RHP}

\subsection{Matrix factorization in terms of the Baker-Akhiezer function}

We are interested in factorizing the matrix $G(t)$ in terms of the function $\CF(z,w)$. In other words, we wish to express two matrices $X^+(t)$ and $X^-(t)$ 
through the Baker-Akhiezer function
such that 
\beq
G(t)=X^+(t)[X^-(t)]^{-1}, \quad t\in L,
\label{5.1}
\eeq
where $X(z)=X^\pm(z)$, $z\in D^\pm$, and $X^+(z)$ and $X^-(z)$ are analytic and nonsingular everywhere
in $D^+$ and $D^-$, respectively, apart from at most a finite number of points where
they may have poles or where $\det X(z)=0$. 
Let $\Gc(z,w)$ be a nontrivial solution to the following homogeneous RHP problem on the surface $\CR$:

 {\it Find a piece-wise meromorphic function $\Gc(z,w)$ with the discontinuity contour $\CL$, $\hat H$-continuous up to the contour
$\CL$ except for a finite number of poles and satisfying the boundary condition
\beq
\Gc^+(t,\xi)=\Gl(t,\xi)\Gc^-(t,\xi), \quad (t,\xi)\in\CL.
\label{5.5}
 \eeq
} 

Then the matrix of factorization $X(z)$ can be expressed exclusively through the function
$\Gc(z,w)$ and the matrix $Y(z,w)$ given by
\beq
Y(z,w)=\fr12\left[I+\fr{1}{h(z)w}Q(z)
\right], \quad Q(z)=\left(\begin{array}{cc}
l_0(z) & l_1(z)\\
l_2(z) & -l_0(z)\end{array}\right), \quad I=\diag\{1,1\},
\label{5.6}
\eeq
in the form [\ref{moi1}], [\ref{ant6}]
\beq
X(z)=\Gc(z,w)Y(z,w)+\Gc(z,-w)Y(z,-w).
\label{5.7}
\eeq
It is a simple matter to verify that
$$
[X(z)]^{-1}=\fr{Y(z,w)}{\Gc(z,w)}+\fr{Y(z,-w)}{\Gc(z,-w)},
$$
\beq
Q^2(z)=h^2(z)f(z)I, \quad Y^2(z,w)=Y(z,w), \quad Y(z,w)Y(z,-w)=0,
\label{5.8}
\eeq
and because of (\ref{5.5}) 
\beq
X^+(t)[X^-(t)]^{-1}=\fr12[\Gl_1(t)+\Gl_2(t)]I+\fr{1}{2h(t)\sqrt{f(t)}}[\Gl_1(t)-\Gl_2(t)]Q(t)=G(t), \quad t\in L.
\label{5.9}
\eeq
We assert that the function  $\Gc_0(z,w)\CF(z,w)$ meets the boundary condition (\ref{5.5}),
and it is bounded at the infinite points $\infty_1$ and $\infty_2$ (the Baker-Akhiezer function $\CF(z,w)$ annihilates the essential singularities 
of the function $\Gc_0(z,w)$  at the infinite points). Thus, the function
$$
\Gc(z,w)=\Gc_0(z,w)\CF(z,w)
$$
\beq
=e^{\Gb_0(z)+w\Gb_1(z)+\GO(P)} \fr{\Gt(u_1(P)-\Gs_1+V_1^\circ, \ldots, u_\Gr(P)-\Gs_\Gr+V_\Gr^\circ)} 
{\Gt(u_1(P)-\Gs_1, \ldots, u_\Gr(P)-\Gs_\Gr)} 
\label{5.9'}
\eeq
is a meromorphic solution to the scalar RHP (\ref{5.5}) on the surface $\CR$, and the matrix (\ref{5.7})
generates  Wiener-Hopf matrix-factors of the matrix $G(t)$.

\subsection{General solution to the scalar RHP on the Riemann surface}

To derive the general solution to the vector RHP (\ref{2.2}), we solve the 
scalar RHP on the Riemann surface $\CR$ (\ref{2.11}). On employing the factorization (\ref{5.5}) of the function $\Gl(t,\Gx)$
we write
\beq
\fr{\Psi^+(t,\Gx)}{\Gc^+(t,\Gx)}=\fr{\Psi^-(t,\Gx)}{\Gc^-(t,\Gx)}+\fr{g^*(t,\Gx)}{\Gc^+(t,\Gx)}, \quad (t,\Gx)\in\CL.
\label{5.10}
\eeq
Since $g^*(t,\Gx)[\Gc^+(t,\Gx)]^{-1}$ is an $\hat H$-continuous function on the surface $\CR$, due to the Sokhotski-Plemelj formulas
it admits a representation in terms of the limit  values of the Weierstrass integral
\beq
F(z,w)=\fr{1}{2\pi i}\int_\CL\fr{g^*(t,\Gx)}{\Gc^+(t,\Gx)}dW, \quad (z,w)\in\CR\setminus\CL,
\label{5.11}
\eeq
as follows:
\beq
F^+(t,\Gx)-F^-(t,\Gx)=\fr{g^*(t,\Gx)}{\Gc^+(t,\Gx)}, \quad (t,\Gx)\in\CL.
\label{5.12}
\eeq
The integral (\ref{5.11}) can be conveniently written as
\beq
F(z,w)=F_1(z)+wF_2(z),
\label{5.13}
\eeq
where
\beq
F_1(z)= \fr{1}{4\pi i}\int_\CL\fr{g^*(t,\Gx)dt}{\Gc^+(t,\Gx)(t-z)}, \quad
F_2(z)= \fr{1}{4\pi i}\int_\CL\fr{g^*(t,\Gx)dt}{\Gx(t)\Gc^+(t,\Gx)(t-z)}.
\label{5.13'}
\eeq
Consequently we may replace the boundary condition (\ref{5.10}) by
\beq
\fr{\Psi^+(t,\Gx)}{\Gc^+(\Gx)}-F^+(t,\Gx)= \fr{\Psi^-(t,\Gx)}{\Gc^-(\Gx)}-F^-(t,\Gx), \quad (t,\Gx)\in\CL,
\label{5.13''}
\eeq
and apply the Liouville theorem
to obtain
\beq
\Psi(z,w)=\Gc(z,w)[F(z,w)+R(z,w)], \quad (z,w)\in{\CR},
\label{5.14}
\eeq
where $R(z,w)$ is a rational function on the surface $\CR$,
\beq
 R(z,w)=R_1(z)+wR_2(z),
\label{5.15}
\eeq
and $R_1(z)$ and $R_2(z)$ are rational functions in the $z$-plane. The function $\Psi(z,w)$ has poles 
at the points with affixes $p_1,p_2,\ldots,p_l$ of   multiplicity  $m_1,m_2,\ldots, m_l$, respectively,  lying in both sheets of the surface. Therefore
 the rational function 
$R(z,w)$ has also poles of the same multiplicity at these points. In addition, due to $\Gr$ simple zeros $Q_j=(t_j,w_j)$ ($j=1,2,\ldots,\Gr$)   of the Baker-Akhiezer function
(these zeros are to be determined) the function $R(z,w)$
may have simple poles  at these points. If $\Gk_m>0$ ($m=1,2$), the function $R(z,w)$ has simple
poles at the points $q_{m j}=(z_{mj}, (-1)^{m-1}\sqrt{f(z_{mj})})\in{\Bbb C}_m\setminus L_m$. Otherwise, if $\Gk_m\le 0$ ($m=1,2$), 
the function $R(z,w)$ is bounded at the points $q_{mj}$. 
Also, the function $R(z,w)$ may have simple poles at the branch points $r_1, r_2,\ldots, r_{2\Gr+2}$, 
the poles of the function $\Psi(z,w)$. Since $|\Psi(z,w)|<c_j|z|^\Gd$ as $z\to\infty_m$, $m=1,2$, $\Gd=\max\{0,\Gd_1-N-\Gr-1\}$,
the functions $R_1(z)$  and $R_2(z)$ may have poles of order $\Gd_1-N-\Gr-1$ and $\Gd_1-N-2\Gr-2$, respectively,  at the infinite point.
The most general form of the rational functions $R_1(z)$  and $R_2(z)$ with the poles
described is given by
$$
R_1(z)=C_0+\sum_{j=1}^{\Gd_1-N-\Gr-1} M_j'z^j+\sum_{j=1}^\Gr\fr{C_jw_j}{z-t_j}
+\sum_{k=1}^l\sum_{j=1}^{m_k}
\fr{D_{kj}'}{(z-p_k)^j}+
\sum_{m=1}^2(-1)^{m-1}\sum_{j=1}^{\Gk_m}\fr{E_{mj}\sqrt{f(z_{mj})}}{z-z_{mj}},
$$
\beq
 R_2(z)=\sum_{j=0}^{\Gd_1-N-2\Gr-2} M_j''z^j+\sum_{j=1}^\Gr\fr{C_j}{z-t_j}+\sum_{j=1}^{2\Gr+2}\fr{K_j}{z-r_j}
+\sum_{k=1}^l\sum_{j=1}^{m_k}
\fr{D_{kj}''}{(z-p_k)^j}+
\sum_{m=1}^2\sum_{j=1}^{\Gk_m}\fr{E_{mj}}{z-z_{mj}},
\label{5.16}
\eeq
where $M_j'$ ($j=1,2,\ldots, \Gd_1-N-\Gr-1)$, $M_j''$ ($j=0,1,\ldots,\Gd_1-N-2\Gr-2$), $C_j$, ($j=0,1,\ldots,\Gr$), $K_j$ ($j=1,2,\ldots,2\Gr+2$), 
$D_{kj}'$, $D_{kj}''$ ($j=1,2,\ldots,m_k$; $k=1,2,\ldots,l$), and $E_{mj}$ ($j=1,2,\ldots,\Gk_m$, $m=1,2$) 
are arbitrary constants. In total, the rational function $R(z,w)$ possesses $\Gk$ free  constants, and $\Gk$ is defined by
$$
\Gk=\left\{\begin{array}{cc}
2\Gd_1+\tilde\Gk+1, \:  &N+2\Gr+2\le  \Gd_1\le 2N+2\Gr+2,\\ 
\Gd_1+N+2\Gr+\tilde\Gk+2, \; & N+\Gr+2\le \Gd_1\le N+2\Gr+1,\\
2N+3\Gr+\tilde\Gk+3, \; & 0\le \Gd_1\le N+\Gr+1,\\
\end{array}
\right.
$$
\beq
 \tilde\Gk=\max\{0,\Gk_1\}+\max\{0,\Gk_2\}.
\label{5.16'}
\eeq

Analysis of formulas (\ref{5.15}) and (\ref{5.16})
shows that the function $R(z,w)$ has simple poles at the points $Q_j=(t_j, w_j)$
and removable singularities at the points $(t_j, -w_j)\in\CR$ ($j=1,2,\ldots,\Gr$). Also, if $\Gk_m>0$,  it has simple poles at the points 
$q_{mj}\in{\Bbb C}_{m}$,  $q_{mj}=(z_{mj}, (-1)^{m-1} \sqrt{f(z_{mj})})$ and removable singularities at the points $(z_{mj}, (-1)^m \sqrt{f(z_{mj})}\in{\Bbb C}_{3-m}$
($j=1,2,\ldots, \Gk_m$;  $m=1,2$).

Owing to the poles of the function $\Gc(z,w)$ and the structure of the functions $F(z,w)$ and $R(z,w)$ we may expect that the function $\Psi(z,w)$
possesses some poles unacceptable  for the solution to the RHP  (\ref{2.11}). Such singularities have to be removed.
Due to the simple poles  of the Baker-Akhiezer function and therefore the poles
of the function $\Gc(z,w)$ at the points $P_1,P_2,\ldots,P_\Gr$ lying
in the first sheet the function $\Psi(z,w)$ has simple poles at these points.  We put
\beq
F(z,w)+R(z,w)=0, \quad (z,w)=(\Gz_j,\sqrt{f(\Gz_j)}), \quad j=1,2,\ldots,\Gr,
\label{5.18}
\eeq
and the points $P_j$ ($j=1,\ldots,\Gr$) become removable singularities.

If $\Gk_m<0$ ($m=1,2$), then the function $\Gc(z,w)$ has $-\Gk_m$ simple poles at the points $q_{mj}\in{\Bbb C}_m$.
For the purpose of removing these poles we request
 \beq
F(z,w)+R(z,w)=0, \quad (z,w)=(z_{mj},(-1)^{m-1}\sqrt{f(z_{mj})}), \quad j=1,2,\ldots,-\Gk_m,\quad m=1,2.
\label{5.19}
\eeq

If $0\le\Gd_1\le N+\Gr+1$, then the function $\Psi(z,w)$ has to be bounded at infinity.
However, due  to the function $w$
in the representations (\ref{5.13}) and (\ref{5.15}) it has order-$\Gr$ poles at the points $\infty_1$ and $\infty_2$. 
Expand the function $F_2(z)+R_2(z)$ in a neighborhood of the infinite point
\beq
F_2(z)+R_2(z)=\fr{v_1}{z}+\ldots+\fr{v_{\Gr}}{z^\Gr}+\fr{v_{\Gr+1}}{z^{\Gr+1}}+\ldots.
\label{5.20}
\eeq
These poles become removable singularities of the function $\Psi(z,w)$ if 
and only if
\beq
v_1=v_2=\ldots=v_\Gr=0.
\label{5.21}  
\eeq
In the case $N+\Gr+2\le\Gd_1\le N+2\Gr$ we have to have  $|\Psi(z,w)|<c_j|z|^{\Gd_1-N-\Gr-1}$ as $z\to\infty_m$, $m=1,2$.
However, the function $\Psi(z,w)$ found has poles of order $\Gr$ at the points $\infty_1$ and $\infty_2$.
Since $1\le\Gd_1-N-\Gr-1\le \Gr-1$, to have the asymptotics required, we have to put 
\beq
v_1=v_2=\ldots=v_{2\Gr+N-\Gd_1+1}=0.
\label{5.22}
\eeq
In the case $N+2\Gr+1\le\Gd_1\le 2N+2\Gr+2$ the function $\Psi(z,w)$ has the asymptotics we need without any extra conditions.

Denote $\hat\Gk=\max\{0,-\Gk_1\}+\max\{0,-\Gk_2\}$. We have $2\Gr+\hat\Gk$, $3\Gr+N-\Gd_1+\hat\Gk+1$ and $\Gr+\hat\Gk$
conditions for the free constants in the cases    $0\le\Gd_1\le N+\Gr+1$, $N+\Gr+2\le\Gd_1\le N+2\Gr$ and  $N+2\Gr+1\le\Gd_1\le 2N+2\Gr+2$, respectively.
If these conditions are fulfilled, then the
function $\Psi(z,w)$ given by (\ref{5.14}) is the general solution to the RHP (\ref{2.11}).

\subsection{Zeros of the Baker-Akhiezer function}

To complete  the procedure presented we have  to  determine the points $Q_j$ ($j=1,\ldots,\Gr$), the zeros of the Baker-Akhiezer
function (\ref{4.30}), or, equivalently, the zeros of the $\Gt$-function (without loss of generality we may assume  that it is not identically equal to zero)
\beq
\Gt(u_1(P)-\Gs_1+V_1^\circ, \ldots, u_\Gr(P)-\Gs_\Gr+V_\Gr^\circ).
\label{6.13.1}
\eeq
We need to know not only the affixes $t_j$ of these zeros, but also to identify the sheet of the surface in which they are located in order to determine
the rational function $R(z,w)$.
On setting
\beq
\Gs_m-V_m^\circ=\sum_{j=1}^\Gr u_m(Q_j)+k_m \quad  ({\rm modulo \;the \;periods}), \quad m=1,2,\ldots,\Gr,
\label{6.14}
\eeq
we obtain that the points $Q_j$ are the zeros of the function $\CF(P)$ indeed.  The system (\ref{6.14}) can equivalently be written as the Jacobi problem of inversion
of abelian integrals:

{\it Find $\Gr$ points on the surface $\CR$, $Q_1,  Q_2,,\ldots, Q_\Gr$, and $2\Gr$ integers, $\mu_1,\mu_2,\ldots, \mu_\Gr$ and $\nu_1,\nu_2,\ldots,\nu_\Gr$, such that
\beq
\sum_{j=1}^\Gr\int_{P_0}^{Q_j}d\Go_m+\sum_{j=1}^\Gr \nu_j B_{mj} +\mu_m=\hat \Gs_m-k_m, \quad m=1,2,\ldots,\Gr,
\label{6.15}
\eeq
where $\hat \Gs_m=\Gs_m-V^\circ_m.$}

This problem reduces [\ref{zve}] to the system of symmetric algebraic equations
\beq
t_1^m+t_2^m+\ldots+t_\Gr^m=\tau_m, \quad m=1,2,\ldots,\Gr,
\label{6.16} 
\eeq
where $\tau_m$ are known and given in terms of the residues at the infinite points [\ref{zve}] or the two zeros of the surface [\ref{ant8}]
of functions expressible in terms of the $\Gt$-function. The system may be converted into the problem of determination
of $\Gr$ zeros of an associated order-$\Gr$ polynomial. 
The integers $\nu_m$ are found by solving the linear system  [\ref{ant8}]
\beq
\sum_{j=1}^\Gr \nu_j \I B_{mj} =\I b _m, \quad m=1,2,\ldots,\Gr,
\label{6.17}
\eeq
while the integers $\mu_m$ are defined by 
\beq
\mu_m=\R b_m-\sum_{j=1}^\Gr \nu_j \R B_{mj}, \quad m=1,2,\ldots,\Gr,
 \label{6.18}
\eeq
explicitly. Here,
\beq
b_m=\hat\Gs_m-k_m-\sum_{j=1}^\Gr u_m(Q_j).
\label{6.19}
\eeq
There are $2^\Gr$ points on the surface $\CR$ which have affixes defined by the $\Gr$ zeros of the polynomial associated with the system (\ref{6.16}).
However, there is one and only one set of points $\{Q_1,\ldots,Q_\Gr\}$ which have the affixes $t_1,\ldots,t_\Gr$, respectively, such that all 
 the numbers $\nu_1,\ldots,\nu_\Gr$ and  $\mu_1,\ldots,\mu_\Gr$  defined by (\ref{6.17}) and (\ref{6.18})  are integers.

\subsection{General solution to the vector RHP}

Having derived the solution to the scalar RHP on the surface $\CR$ (\ref{2.11}) we can now determine and examine the solution to the original
vector RHP (\ref{2.2}). From (\ref{2.7}) and (\ref{2.6}) we express the components of the vector $\BGF(z)$, $\GF_1(z)$ and $\GF_2(z)$,
as
$$
\GF_1(z)=\psi_1(z)+\psi_2(z),
$$
\beq
\GF_2(z)=-\fr{l_0(z)}{l_1(z)}[\psi_1(z)+\psi_2(z)]+\fr{h(z)\sqrt{f(z)}}{l_1(z)}[\psi_1(z)-\psi_2(z)],\quad z\in{\Bbb C},
\label{6.20}
\eeq
where
$\psi_m(z)=\Psi(z,w)$, $(z,w)\in{\Bbb C}_m$, $m=1,2$. 
We have obtained the solution of the RHP (\ref{2.11})  in the class of functions having the poles $p_1, \ldots, p_l$ of multiplicity
$m_1,\ldots,m_l$, respectively, due to the presence of the polynomial $h(z)$ in (\ref{2.7''}) and its zeros at these points. 
However, the solution to the original RHP (\ref{2.2}), the vector $\BGF(z)$,
has to be analytic at these points.   This can be achieved by introducing the following $N$ conditions
\beq
\lim_{z\to p_j}\fr{d^k}{dz^k}\left\{(z-p_j)^{m_j-k}[\psi_1(z)+\psi_2(z)]\right\}=0, \quad k=0,1,\ldots,m_j-1, \quad j=1,2,\ldots,l.
\label{6.22} 
\eeq
If these conditions are satisfied, then the functions $\GF_1(z)$ and $\GF_2(z)$ are analytic at the poles of the functions $\psi_1(z)$
and $\psi_2(z)$ (the zeros of $h(z)$). 

Let $\tilde p_j$ be  order-$\tilde m_j$ zeros  ($j=1,2,\ldots, \tilde l$) of the function $l_1(z)$, $\tilde m_1+\ldots+\tilde m_{\tilde l}=\Gd_1$.
These zeros are poles of the same multiplicity of the function $\GF_2(z)$ in (\ref{6.20}). To remove these poles we require
 \beq
\lim_{z\to\tilde p_j}\fr{d^k}{dz^k}\left[(z-\tilde p_j)^{\tilde m_j-k} \GF_2(z)\right]=0, \quad k=0,1,\ldots,\tilde m_j-1, \quad j=1,2,\ldots,\tilde l.
\label{6.23} 
\eeq

Finally, we need to guarantee that the functions $\GF_1(z)$ and $\GF_2(z)$ are bounded at infinity. 
Analyze first the case $0\le \Gd_1\le N+\Gr$. Since $|\psi_j(z)|\le c_j$ ($j=1,2$) as $z\to\infty$, it follows from (\ref{6.20}) that the function $\GF_1(z)$
is bounded.
 Expand the function $\GF_2(z)$ in a neighborhood of the infinite point
\beq
\GF_2(z)=\tilde v_{N+\Gr-\Gd_1+1}z^{N+\Gr-\Gd_1+1}+\ldots+\tilde v_1z+\tilde v_0+\ldots.
\label{6.24}
\eeq
On putting 
\beq
\tilde v_1=\tilde v_2=\ldots=\tilde v_{N+\Gr-\Gd_1+1}.
 \label{6.25}
 \eeq
we remove the growth  of the function $\GF_2(z)$. As $\Gd_1=N+\Gr+1$, the function
$\GF_2(z)$ is bounded unconditionally.

 Consider now the case  $N+\Gr+2\le \Gd_1\le 2N+2\Gr+2$. It follows from (\ref{6.20}) and the asymptotics of the
 functions $\psi_1(z)$ and $\psi_2(z)$  that at the infinite point the function $\GF_2(z)$ is bounded, while 
 the function $\GF_1(z)$ has a pole of order $\Gd_1-N-\Gr-1$.
  Let
\beq
\psi_1(z)+\psi_2(z)=\hat v_{\Gd_1-N-\Gr-1}z^{\Gd_1-N-\Gr-1}+\ldots+\hat v_1 z+\hat v+\ldots, \quad z\to\infty.
\label{6.26}
\eeq
On satisfying the conditions
\beq
\hat v_1=\hat v_2=\ldots=\hat v_{\Gd_1-N-\Gr-1}=0
 \label{6.27}
 \eeq
we obtain the function $\GF_1(z)$ bounded at the infinite point.

We now summarize the results.

{\sl Theorem. Let $G(t)$ be a nonsingular $2\times 2$ matrix 
\beq
G(t)=\left(\begin{array}{cc}
\Ga_1(t)+\Ga_2(t)l_0(t) & \Ga_2(t)l_1(t) \\
 \Ga_2(t)l_2(t) & \Ga_1(t)-\Ga_2(t)l_0(t)\\
\end{array}\right),
\label{6.27}
\eeq
where 
$\Ga_1(t), \Ga_2(t)l_j(t)\in \hat H(L)$, $j=0,1,2$, $l_0(t)$, $l_1(t)$ and $l_2(t)$ are polynomials, and $L$ is the real axis.
Denote $l_0^2(z)+l_1(z) l_2(z)=h^2(z)f(z)$, $\Gd_1=\deg l_1(z)$  and  $N=\deg h(z)$.
Assume that $2\Gr+2=\deg f(z)$, the zeros of the polynomial $f(z)$ are simple,  and none of the zeros of $f(z)$ and $h(z)$
fall in the contour $L$. 

Let $\Gk_1$ and $\Gk_2$ be the integers defined by
$\Gk_j=\ind\Gl_j(t)$, $t\in L$, where $\Gl_1(t)$ and $\Gl_2(t)$ are the eigenvalues of $G(t)$,   $\Gl_j=\Ga_1-(-1)^j \Ga_2 h\sqrt{f}$, $j=1,2$.
Denote $ \tilde\Gk=\max\{0,\Gk_1\}+\max\{0,\Gk_2\}.$

Then the functions (\ref{6.20}) possess 
$\Gk$ arbitrary constants
\beq
\Gk=\left\{\begin{array}{cc}
2\Gd_1+\tilde\Gk+1, \:  &N+2\Gr+2\le  \Gd_1\le 2N+2\Gr+2,\\ 
\Gd_1+N+2\Gr+\tilde\Gk+2, \; & N+\Gr+2\le \Gd_1\le N+2\Gr+1,\\
2N+3\Gr+\tilde\Gk+3, \; & 0\le \Gd_1\le N+\Gr+1,\\
\end{array}
\right.
\label{6.28}
\eeq
 which have to satisfy $\Gk'=\Gk-\Gk_1-\Gk_2-2$
additional conditions  (\ref{5.18}),  (\ref{5.19}), (\ref{5.21}),  (\ref{5.22}),
(\ref{6.22}), (\ref{6.23}), (\ref{6.25}) and (\ref{6.27}). If $\Gk_1+\Gk_2\ge -2$, then
 the solution to the problem (\ref{2.2}) exists, has $\Gk_1+\Gk_2+2$ free constants and is defined by (\ref{6.20}).
Otherwise, the solution does not exist. If however the 
vector $\Bg(t)$  satisfies $-\Gk_1-\Gk_2-2$ conditions which guarantee that all the additional conditions  are fulfilled, then the solution exists  and it is unique.}

\section*{Conclusions}

We have proposed a new technique for deriving  Wiener-Hopf factors  of the Chebotarev-Khrapkov matrix 
$G(t)=\Ga_1(t)I+\Ga_2(t)Q(t)$, $\Ga_1(t), \Ga_2(t) Q(t)\in \hat H(L)$, $Q(t)$ is a $2\times 2$ zero-trace polynomial matrix. The method has been applied to solve
the vector RHP $\BGF^+(t)=G(t)\BGF^-(t)+\Bg(t)$, $t\in L$. The known technique [\ref{moi1}], [\ref{ant6}] first
reduces the vector problem to a scalar RHP on the Riemann surface $\CR$ of the algebraic function $w^2=f(z)$, $\det Q(z)=h^2(z)f(z)$.  Then it 
 finds a function $\Gc_0(z,w)$ which factorizes
the coefficient of the RHP on the surface and allows for  essential singularities at the infinite points of $\CR$. These singularities are removed
by solving a certain Jacobi problem of inversion of hyperelliptic integrals. At this stage, a meromorphic solution is derived.   
The inadmissible poles due to the technique applied are removed afterwards.
In contrast with this method, the technique we have developed hinges on the derivation of the Baker-Akhiezer function widely used
in the theory of integrable systems. This procedure quenches the essential singularities by constructing a special abelian integral
of the second type $\GO(P)$. It has zero $A$-periods, and the  principal part of the function $\exp\{\GO(P)\}$ at the infinite points 
is derived according to the behavior of
the function $\Gc_0(z,w)$ at the infinite points. The consequent use of the quotient of two Riemann $\Gt$-functions
serves to annihilate the discontinuity of the function $\exp\{\GO(P)\}$ due to the nonzero $B$-periods of the integral $\GO(P)$.
The product of the function $\exp\{\GO(P)\}$ and the quotient of the two $\Gt$-functions forms the Baker-Akhiezer function $\CF(P)$, while
the product of $\Gc_0(z,w)$ and  $\CF(P)$ forms a solution of the homogeneous scalar RHP on the surface. It does not have essential
singularities and is a meromorphic function in $\CR$ with a finite number of prescribed poles. This gives Wiener-Hopf factors of $G(t)$
and does not require the solution of a Jacobi inversion problem. For the general solution of the vector RHP however the solution of the associated 
Jacobi inversion problem is unavoidable. This is because the Baker-Akhiezer function has $\Gr$ zeros ($\Gr$ is the genus of the surface $\CR$),
and their location cannot be prescribed. At the stage of application of the generalized Liouville theorem, the zeros of the Baker-Akhiezer function are 
 needed for determination of the rational vector in the general solution. This information can be recovered by
stating and solving the corresponding Jacobi inversion problem.

 \vspace{.2in}

{\large {\bf  Acknowledgments.}} 
The author is thankful to A. R. Its for discussions of the paper.

 \vspace{.2in}

{\centerline{\Large\bf  References}}

\begin{enumerate}

\item\label{ant1}
Antipov, Y. A.: An exact solution of the 3-D-problem on an interface semi-infinite plane crack.  J.
Mech. Phys. Solids 47 1051-1093 (1999) 

\item\label{ant2}
Antipov, Y. A.: Solution by quadratures of the problem of a cylindrical crack by the method of matrix factorization. IMA J. Appl. Math.
66 591-619 (2001)

\item\label{ant3}
Antipov,  Y. A.: A symmetric Riemann-Hilbert problem for order-4 vectors in diffraction theory.
Quart. J. Mech. Appl. Math. 63 349-374 (2010)

\item\label{ant4}
 Antipov,  Y. A.: A genus-3 Riemann-Hilbert problem and diffraction of a wave by orthogonal resistive half-planes.
 Comput. Meth.  Function Theory 11 439-462 (2011)

\item\label{ant5}
Antipov, Y. A. and  Moiseev, N. G.:  Exact solution of the plane problem for a composite plane with a cut
across  the boundary between two media. J. Appl. Math. Mech. 55 531-539 (1991)

\item\label{ant6}
 Antipov, Y. A. and  Silvestrov, V. V.:  Factorization on a Riemann surface in scattering theory. 
  Quart. J. Mech. Appl. Math. 55 607-654 (2002)

\item\label{ant7}
 Antipov, Y. A. and  Silvestrov, V. V.: Vector functional-difference equation in electromagnetic scattering.
   IMA J. Appl. Math. 69 27-69 (2004)

\item\label{ant8} 
Antipov, Y. A. and   Silvestrov, V. V.:
Electromagnetic scattering from an anisotropic impedance half plane
at oblique incidence: the exact solution.  Quart. J. Mech. Appl. Math. 59 211-251 (2006)

\item\label{buy}
B\"uy\"ukaksoy, A. and Serbest, A. H.:  Matrix Wiener-Hopf methods applications to some diffraction problems. In: 
Hashimoto, M. Ideman, M.,  Tretyakov, O. A.(eds), 
 Analytical and numerical methods in electromagnetic wave theory. 257-315, Science House Co. Ltd, Tokyo (1993)

\item\label{che1}
Chebotarev,  G. N.:  On closed-form solution of a Riemann boundary value problem for n pairs
of functions.  Uchen. Zap. Kazan. Univ. 116 31-58 (1956)

\item\label{che2}
Chebotarev, N. G.: Theory of algebraic functions.  OGIZ, Moscow (1948)

\item\label{dan}
Daniele, V. G.: On the solution of vector Wiener-Hopf equations occurring in scattering problems.  Radio Sci.
 19 1173-1178 (1984)

\item\label{dub1}
Dubrovin, B. A.: The inverse scattering problem for periodic finite-zone
potentials. Funct. Anal. Appl. 9 61-62 (1975)

\item\label{dub2}
Dubrovin, B. A.: Theta functions and non-linear equations.  Russian Math. Surveys 36:2 11-92 (1981)

\item\label{dub3}
Dubrovin, B. A.,  Matveev, V. B. and Novikov, S. P.: Nonlinear equations of Korteweg-de Vries type, finite-band linear operators and Abelian varieties.
Russian Math. Surveys 31:1 59-146 (1976)

\item\label{hur}
Hurd, R. A. and   L\"uneburg, E.: Diffraction by an anisotropic impedance half plane.  Canad. J.
Phys.  63 1135-1140 (1985)

\item\label{its}
Its, A. R. and Matveev, V. B.: Schr\"odinger operators with the finite-band spectrum and the $N$-soliton solutions of the Korteweg - de Vries  equation. 
 Theor. Math. Phys. 23 343-355 (1975)

\item\label{jon}
Jones, C. M. A.: Scattering by a semi-infinite sandwich panel perforated on one side. 
Proc. R. Soc. A 454 465-479 (1990)

\item\label{khr}
Khrapkov, A. A.: Certain cases of the elastic equilibrium of an infinite wedge with a nonsymmetric 
 notch at the vertex, subjected to concentrated forces.
 J. Appl. Math. Mech. 35 625-637 (1971)

\item\label{kri}
Krichever, I. M.:  Methods of algebraic geometry in the theory of nonlinear equations.  Russian Math. Surveys 32:6 185-213 (1971)

\item\label{lun} 
L\"uneburg, E. and Serbest, A. H.: Diffraction of an obliquely incident plane wave by a two-face
impedance half plane: Wiener-Hopf approach.  Radio Sci.  35 1361-1374 (2000)

\item\label{mat}
Matveev, V. B.:  30 years of finite-gap integration theory.  Phil.
Trans. R. Soc. A  366 837-875 (2008)

\item\label{moi1}
Moiseev, N. G.:  Factorization of matrix functions of special form. Soviet Math. Dokl.
39 264-267 (1989)

\item\label{moi2}
Moiseyev, N. G. and  Popov, G. Ya.: Exact solution of the problem of bending of a semi-infinite
plate completely bonded to an elastic half-space.
Izv. Akad. Nauk SSSR, Solid Mechanics 25 113-125 (1990)

\item\label{raw}
Rawlins, A.D.: Two waveguide trifurcation problems. Math. Proc. Camb. Phil. Soc.
121 555-573 (1995)

\item\label{spr} 
Springer, G.:  Introduction to Riemann Surfaces.  Addison-Wesley, Reading (1956) 

\item\label{vek}
Vekua, N. P.: Systems of Singular Integral Equations. Noordhoff, Groningen (1967)

\item\label{zve}
Zverovich, E. I.: Boundary value problems in the theory of analytic functions in H\"older classes
on Riemann surfaces. Russian Math. Surveys 26:1 117-192 (1971)

\item\label{zve2}
Zverovich, E. I.: The problem of linear conjugation on a closed Riemann surface.
Compl. Anal. Oper. Theory  2  709-732 (2008)

\end{enumerate}

\end{document}